\magnification=\magstep1
\overfullrule=0pt
\font\caps=cmcsc10
\def\complex{\mathop{{\rm I}\kern-.6em{\bf C}}\nolimits}
\def\nat{\mathop{{\rm I}\kern-.2em{\rm N}}\nolimits}
\def\que{\mathop{{\rm I}\kern-.6em{\rm Q}}\nolimits}
\def\real{\mathop{{\rm I}\kern-.2em{\rm R}}\nolimits}
\def\D{{\cal D}} 
\def\bb{{\bf b}} 
\def\bff{{\bf f}} 
\def\bg{{\bf g}} 

\def\dis{\mathop{\rm dis}\nolimits}
\def\Im{\mathop{\rm Im}\nolimits}
\def\osc{\mathop{\rm osc}\nolimits}
\def\Re{\mathop{\rm Re}\nolimits}
\def\supp{\mathop{\rm supp}\nolimits}
\def\varep{\varepsilon}
\def\frac#1#2{{\textstyle{#1\over#2}}} 
\def\disteq{\mathop{\buildrel {\hbox{\sevenrm dist}} \over =}\nolimits}
\def\subsetneqq{\ \lower.5ex\hbox{${\buildrel{\scriptstyle\subset}\over
	{\scriptstyle\ne}}$}\ }
\def\blackbox{\hbox{\vrule width6pt height7pt depth1pt}}
\def\qed{~\hfill~\blackbox\medskip}		
\def\myskip{\noalign{\vskip6pt}}	
\def\newpage{\vfill\eject} 
\def\iitem{\itemitem}
\def\hangbox to #1 #2{\vskip1pt\hangindent #1\noindent \hbox to #1{#2}$\!\!$}
\def\ref#1{\hangbox to 30pt {#1\hfill}}
\outer\def\beginsection#1\par{\vskip0pt plus.3\vsize
	\vskip0pt plus-.3\vsize\bigskip\vskip\parskip
	\message{#1}\leftline{\bf#1}\nobreak\smallskip\noindent} 
\topinsert\vskip.5in\endinsert
\centerline{\bf On Weakly Null FDD's in Banach Spaces}
\medskip
\centerline{\sl by}
\medskip
\centerline{{\caps E. Odell,}\footnote*{Research partially supported by 
the National Science Foundation and TARP 235.}
{\caps H.P. Rosenthal$^*$ and Th. Schlumprecht}}
\vskip.2in

{\narrower\smallskip\noindent 
{\bf Abstract.}
In this paper we show that every sequence $(F_n)$ of finite dimensional 
subspaces of a real or complex Banach space with increasing dimensions can be 
``refined'' to yield an F.D.D.\ $(G_n)$, still having increasing dimensions, 
so that either every bounded sequence $(x_n)$, with $x_n \in G_n$ for 
$n\in\nat$, is weakly null, or every normalized sequence $(x_n)$, with 
$x_n \in G_n$ for $n\in\nat$, is equivalent to the unit vector basis 
of $\ell_1$. 

Crucial to the proof are two stabilization results concerning Lipschitz 
functions on finite dimensional normed spaces. These results also lead 
to other applications. We show, for example, that every infinite 
dimensional Banach space $X$ contains an F.D.D.\ $(F_n)$, with $\lim_{n\to
\infty} \dim (F_n)=\infty$, so that all normalized sequences $(x_n)$, with 
$x_n\in F_n$, $n\in\nat$, have the same spreading model over $X$. 
This spreading model must necessarily be 1-unconditional over $X$. 
\smallskip}

\baselineskip=18pt		
\beginsection{\S1. Introduction} 

Let $(F_n)$ and $(G_n)$ be two sequences of finite dimensional subspaces of 
a Banach space $X$. We say $(F_n)$ is {\it large\/} if $\lim_{n\to\infty} 
\dim F_n =\infty$. We say $(G_n)$ is a {\it refinement\/} of $(F_n)$ 
if there is a strictly increasing sequence $(k_n)\subset \nat$ so that 
$G_n$ is a subspace of $F_{k_n}$ for all $n\in \nat$. If each $(F_n)$ has 
a given basis $\bb_n = (f_i^{(n)} : 1\le i\le \dim F_n)$, we say $(G_n)$ 
is a {\it block refinement\/} of $(F_n)$ with respect to $(\bb_n)$ if 
$G_n$ is spanned by a block basis of $\bb_n$ for all $n$. 
$(F_n)$ is called an F.D.D. (Finite Dimensional Decomposition) if $(F_n)$ 
is a Schauder-decomposition for its closed linear span. It is readily seen 
(using the standard Mazur argument) that every large sequence $(F_n)$ 
has a large F.D.D.\ refinement $(G_n)$; moreover $(G_n)$ can be chosen 
to be a block-refinement of $(F_n)$ with  respect to $(\bb_n)$ for a given 
sequence of bases $(\bb_n)$ of the F.D.D. 
We say $(G_n)$ is {\it weakly null\/} if every bounded sequence $(x_n)$ with 
$(x_n) \in G_n$ for all $n$, is weakly null. We say $(G_n)$ is 
{\it uniformly\/}-$\ell_1$ if there exists a $C>0$ such that all normalized 
sequences $(x_n)$ with $x_n \in G_n$ for all $n$, are $C$-equivalent  to 
the unit vector basis of $\ell_1$. Of course $(G_n)$ is uniformly-$\ell_1$ 
precisely when $(G_n)$ is an $\ell_1$-F.D.D.; that is, the closed linear 
span of the $G_n$'s is canonically isomorphic to $(\sum \oplus G_n)_1$, the 
space of all sequences $(g_n)$ with $g_n\in G_n$ for all $n$ and $\|(g_n)\| 
\buildrel {\rm df}\over = \sum \|g_n\| <\infty$. 

Except as noted, our terminology is standard and may be found in the 
book [LT]. All Banach spaces are assumed to be separable. 

If $(x_n)$ (resp.\ $(G_n)$) is a (finite or infinite) sequence of elements 
of (resp.\ finite-dimensional subspaces of) a Banach space $X$, 
$[x_n]$ (resp.\ $[G_n]$) denotes the closed linear span of $(x_n)$ 
(resp.\ $(G_n)$). $S_X$ denotes the unit sphere of $X$ and $Ba (X)$ its 
unit ball. 

Our main result is the following. 

\proclaim Theorem 1. 
Let $(F_n)$ be a large sequence of finite-dimensional subspaces of a Banach 
space $X$. Then there exists a large refinement $(G_n)$ of $(F_n)$ so 
that either $(G_n)$ is a weakly null {\rm FDD} or $(G_n)$ is an 
$\ell_1$-{\rm FDD}. Furthermore if there is a given sequence $(\bb_n)$ 
of bases of the $F_n$'s with uniformly bounded basis constants, then the 
above sequence $(G_n)$ can be chosen to be a block refinement of $(F_n)$ 
with respect to $(\bb_n)$. 

Theorem 1 can be viewed as a block version of the $\ell_1$-theorem of the 
second named author, which says that every normalized sequence $(x_n)$ 
in a Banach space $X$ has a subsequence which is either equivalent to the 
unit vector basis of $\ell_1$ or is weak Cauchy [R1]. Using Krivine's 
theorem [K] (which is also used the in proof of Theorem~1), one gets 
further structural consequences of this block version. 
Krivine's theorem (as refined in [R2] and finally in [L]) may be 
formulated as follows: 
\smallskip

{\narrower\smallskip\noindent \sl 
Given a large sequence $(F_n)$ of finite-dimensional subspaces of a Banach 
space with bases $(\bff_n)$ with uniformly bounded basis constants, there 
exists a block refinement $(G_n)$ of $(F_n)$ with block bases $(\bg_n)$ of 
the $\bff_n$'s so that for all $n$, $n=\dim (G_n)$ and $\bg_n$ is 
$1+{1\over n}$-equivalent to the unit vector basis of $\ell_p^n$. 
\smallskip}

\noindent  Of course it thus follows that the $G_n$'s in the conclusion of 
Theorem~1 can be chosen to be uniformly isomorphic to $\ell_p^n$, for some 
$1\le p\le\infty$. We thus obtain immediately the following  result. 

\proclaim Corollary 2. 
Let $(F_n)$ be a large sequence of finite dimensional subspaces of a 
Banach space $X$, with given bases $(\bb_n)$ with uniformly bounded basic 
constants; and assume no normalized sequence $(f_n)$ with $f_n\in F_n$ for 
all $n$, has a weak Cauchy subsequence. Then there exists $1\le p\le\infty$ 
and a block refinement $(G_n)$ of $(F_n)$ with respect to $(\bb_n)$, such 
that $[G_n]$ is canonically isomorphic to $(\sum \oplus \ell_p^n)_1$. 

Now Corollary 2 trivially implies that {\sl if  $X$ has the Schur property 
and contains $\ell_p^n$'s uniformly, then $(\oplus \ell_p^n)_1$ 
embeds in\/} $X$. 
Of course this is trivial if $1\le p\le 2$, since then $\ell_p$ is finitely 
represented in $\ell_1$. However the following immediate block version 
does not appear to be obvious for any value of $p$ larger than $1$. 

\proclaim Corollary 3. 
Let $X$ have the Schur property, and suppose, for some $1<p\le\infty$, that 
$\ell_p$ is block finitely represented in a particular basic sequence  
$(x_j)$ in $X$. Then some block basis of $(x_j)$ is equivalent to the 
natural basis of $(\sum \oplus \ell_p^n)_1$. 

A famous question in Banach space theory was whether any infinite dimensional 
Banach space $X$ which does not contain $\ell_1$ isomorphically  must contain 
an infinite-dimensional subspace with a separable dual. This is equivalent 
to asking whether such an $X$ contains a shrinking basic sequence $(x_n)$; 
i.e., a basic sequence $(x_n)$ so that each bounded block basis $(y_n)$ 
is weakly null. Of course if $(x_n)$ is such a sequence and $(k_n)$ is 
an increasing sequence in $\nat\cup \{0\}$ with $k_{n+1}-k_n\to\infty$, 
then setting $F_n = [x_i]_{i=k_n+1}^{k_{n+1}}$, $(F_n)$ is a large weakly 
null FDD. However T.~Gowers [G2] has recently solved the general problem 
in the negative; i.e., there is a Banach space $X$ not containing $\ell_1$, 
with {\it no\/} shrinking basic sequences. Nevertheless, Theorem~1 gives at 
once  that {\sl every basic sequence in any $X$ not containing $\ell_1$ 
has a block basis $(x_n)$ which yields large weakly null FDD's 
as above\/}. 

\proclaim Corollary 4. 
If $\ell_1$ is not isomorphically contained in $X$ and $(y_n)$ is a 
basic sequence in $X$, then for each increasing sequence $(k_n)\subset 
\nat \cup\{0\}$ there exists a block basis $(x_n)$ of $(y_n)$ so that 
$(F_n)$ is weakly null, where $F_n = [x_i]_{i=k_n+1}^{k_{n+1}}$ for 
all $n$. 

Corollary 4 motivates the following problem. 
\medskip

\noindent {\bf Problem.} 
Assume $\ell_1$ is not contained in $X$. Does there exist a basic sequence 
$(x_n)$ so that all bounded ``admissible'' block bases of $(x_n)$ 
converge weakly to zero? (We call a block basis $(y_n)$ of $(x_n)$ 
{\it admissible\/} if $y_n = \sum_{i=1}^{\ell_n} \alpha_i^{(n)} x_{m_i^{(n)}}$ 
where  for each $n$, $\ell_n\in\nat$, $(\alpha_i^{(n)})\in\real^{\ell_n}$ and 
$\ell_n \le m_1^{(n)} < m_2^{(n)} <\cdots < m_{\ell_n}^{(n)}$. 
In the terminology of [FJ] this just says that for all $n$, 
$\supp (y_n)$ (with respect to $(x_n)$) is an admissible subset of $\nat$.) 
\medskip

Another corollary of Theorem 1 is the following result, stated in [R6,  
Corollary~22] and proved there using Theorem~1 and the Borsuk antipodal 
mapping theorem. Corollary~5 was obtained independently by 
W.B.~Johnson and T.~Gamelin [CGJ]. 

\proclaim Corollary 5. 
Assume $\ell_1$ does not embed in $X$, where $X$ is an infinite 
dimensional Banach space. 
Then there exists a normalized weakly null basic sequence $(x_i)$ in $X$ 
possessing a normalized sequence of biorthogonal functionals.  

The main tools needed to prove Theorem 1 will be the following two 
finite  dimensional ``stabilization principles.'' The first one was 
observed by V.~Milman (see [MS, p.6]) in connection with A.~Dvoretzky's 
famous theorem that in every infinite dimensional Banach space 
one finds, for each $\varep >0$ and $n\in\nat$, an $n$-dimensional subspace 
$F$ which is $(1+\varep)$-isomorphic to $\ell_2^n$. The second 
stabilization principle follows mainly from Lemberg's [L] proof of 
Krivine's theorem. 
\medskip

\proclaim  First Stabilization Principle. 
\vskip1pt
For every $C>0$, $\varep >0$ and $k\in\nat$ there is an $n=n(C,\varep,k)
\in\nat$ so that: If $F$ is an $n$-dimensional normed space and $f:F\to\real$ 
is $C$-Lipschitz (i.e., $|f(x)-f(y)| \le C\|x-y\|$ for $x,y\in F$), then 
there is a $k$-dimensional subspace $G$ of $F$ so that 
$$\osc \bigl(f|_{S_G}\bigr)\equiv\sup \bigl\{ |f(x)-f(y)| :x,y\in S_G\bigr\}
<\varep\ .$$ 

\proclaim Second Stabilization Principle.
\vskip1pt
For all $C>0$, $\varep>0$ and $k\in\nat$ there is an $n=n(C,\varep,k)\in\nat$ 
so that if $F$ is an $n$-dimensional normed space with a basis $(x_i)_{i=1}^n$, 
whose basis constant does not exceed $C$, and if $f:F\to\real$ is 
$C$-Lipschitz, then there is a block basis $(y_i)_{i=1}^k$ of 
$(x_i)_{i=1}^n$ so that 
$$\osc \Bigl(f|_{S_{[y_i]_{i=1}^k}}\Bigr) <\varep\ .$$ 

Since on the one hand the second stabilization principle nearly follows 
in a straightforward manner from the proof of Krivine's theorem (the only 
exception is the case $F=\ell_\infty^n$), but on the other hand does not 
follow from the statement of Krivine's theorem itself, we will 
sketch the proof in section~3. 

The next result gives another application of the above stabilization 
principles. The result yields that for a given Lipschitz function $f$ 
and large sequence  $(F_n)$ of $X$ of finite-dimensional subspaces, 
there exists a large refinement $(G_n)$, a Banach space $E$ with a 
one-unconditional basis $(e_j)$, and a function $\tilde f :E\to \real$ 
so that for all sequences $(x_i)$ with $f_i \in S_{G_i}$ for all $i$, and all 
$k$, and all sequences $(\alpha_i) \in Ba (\ell_\infty)$ 
$$\tilde f \biggl( \sum_{i=1}^k \alpha_i e_i \biggr) 
= \lim_{n_k >\cdots > n_1\to\infty} 
f\biggl( \sum_{i=1}^k \alpha_i x_{n_i} \biggr)\ .$$ 
The result may be formulated quantitatively as follows: 
($c_{00}$ denotes the linear space of finitely 
supported real valued functions on $\nat$. We write for $A,B\in\real$ 
and $\varep>0$, $A\buildrel \varep\over = B$ if $|A-B|<\varep$.) 

\proclaim Theorem 6. 
Let $X$ be an infinite dimensional Banach space and let $f:X\to \real$ 
be Lipschitz. Let $(\varep_n)\subset \real_+$ with $\lim_{n\to\infty} 
\varep_n=0$ and let $(F_n)$ be a large sequence of finite dimensional subspaces 
of $X$. There exists a large refinement $(G_n)$ of $(F_n)$ and a function 
$\tilde f :c_{00} \cap Ba(\ell_\infty) \to \real$ so that: 
For all $k\in\nat$ and $n_1,n_2,\ldots,n_k\in \nat$ with $k\le n_1< n_2<
\cdots <n_k$, and all $(\alpha_i)_{i=1}^k \in Ba(\ell_\infty^k)$, 
$$\tilde f (\alpha_1,\alpha_2,\ldots,\alpha_k,0,0,\ldots) 
\buildrel {\varep_k}\over = f\biggl( \sum_{i=1}^k \alpha_i x_i\biggr)$$ 
whenever $x_i \in S_{G_{n_i}}$ for $1\le i\le k$.  
Moreover if each $F_n$ has a given basis $\bb_n$ whose basis constant does not 
exceed some fixed number, $(G_n)$ may be chosen to be a block 
refinement of $(F_n)$  with respect to $(\bb_n)$. 

Theorem 6 has a consequence concerning spreading models, and in fact the Banach 
space ``$E$'' given in the above qualitative formulation may be chosen 
to be a spreading model of $X$. 
Recall that (see e.g., [BL] or [O]) every 
seminormalized basic sequence in $X$ admits  a 
subsequence $(x_n)$ satisfying: For all $x\in X$, $k\in\nat$ and 
$(\alpha_i)_{i=1}^k \subseteq \real$, 
$$\lim_{n_1\to\infty}\ \lim_{n_2\to\infty}\ldots \lim_{n_k\to\infty} 
\Big\| x + \sum_{i=1}^k a_i x_{n_i}\Big\| \ \hbox{ exists.}$$ 
The limit is denoted by $\|x+\sum_{i=1}^k a_i e_i\|$ and defines 
a norm on $X\oplus E$ where $E= [e_i]$. $E$ is called a 
{\it spreading model\/} of $X$ and $X\oplus E$ is called a 
{\it spreading model\/} of 
$(x_i)$ over $X$. $(e_i)$ is {\it $1$-unconditional 
over\/} $X$ if for all $x\in X$, $(\alpha_i)_1^k \subseteq\real$ and 
$(\varep_i)$  with $|\varep_i| = 1$ for all $i$, 
$$\Big\| x+\sum_{i=1}^k \alpha_i e_i\Big\| 
= \Big\| x+\sum_{i=1}^k \varep_i \alpha_i e_i\Big\|\ .$$

\proclaim  Corollary 7.
Every large sequence $(F_n)$ of finite dimensional subspaces of an infinite 
dimensional Banach space $X$ 
has a large refinement $(G_n)$ with the following property: All sequences 
$(x_n)$, with $x_n\in S_{G_n}$ for $n\in\nat$, have the same spreading 
model $E= [e_i]$ over $X$. In particular $(e_i)$ is $1$-unconditional over $X$. 
Moreover, $G_n$ can be chosen, so that for all $\varep>0$, $k\in \nat$ and 
$x\in X$ there exists $k_0 \in \nat$ such that if $k_0\le n_1<n_2<\cdots 
< n_k$, then 
$$\left| \Big\| x+\sum_{i=1}^k \alpha_i e_i\Big\| 
- \Big\| x+ \sum_{i=1}^k \alpha_i x_i\Big\| \right| < \varep$$ 
whenever $x_i\in S_{G_{n_i}}$, $i=1,\ldots,k$, and 
$(\alpha_i)_{i=1}^k \in Ba (\ell_\infty^k)$. 

As usual, there is a corresponding ``block refinement'' version. 
Corollary 7 follows from Theorem~6 and a standard diagonal argument 
using the Lipschitz functions $f_x(y) = \|x+y\|$ as $x$ ranges over a 
dense subset of $X$. The result that every Banach space $X$ has a spreading 
model which is 1-unconditional over $X$ is due to the second named 
author, see [R4], [R5]. 

We note finally  an application of Theorems 1 and 6 to the Banach-Saks 
property. The following principle was discovered in 1975 (cf.\ [R2]; 
a proof may be found in [BL]). 
\smallskip

{\narrower\smallskip\noindent 
{\sl Given $(x_j)$ a semi-normalized weakly  null sequence in 
a Banach space, there is a subsequence $(x'_j)$ so that either $(x'_j)$ 
has a spreading model isomorphic to $\ell_1$, or 
$\frac1n \|\sum_{j=1}^n x''_j\| \to0$ as $n\to\infty$ for 
all further subsequences $(x''_j)$ of $(x'_j)$.} 
\smallskip}

Now in fact one may assume in any case that $(x'_j)$ generates a 
spreading model, with basis $(b_j)$ say; then the second alternative occurs 
precisely when $(b_j)$ itself is weakly null. In this case, one has $\|\frac1n 
\sum_{j=1}^n b_j\|\to 0$ as $n\to\infty$. Then, e.g., 
setting $\varep_n = \frac2n \|\sum_{j=1}^n b_j\|$, $(x'_j)$ can be chosen 
so that $\frac1n \|\sum_{j=1}^n x''_j\|\le \varep_n$ for all subsequences 
$(x''_j)$ of $(x'_j)$. 

The following result now follows immediately from Theorem~1 and Corollary~7. 

\proclaim Corollary 8. 
Let $(F_j)$ be a large sequence of finite dimensional subspaces of a  
Banach space $X$, so that no normalized sequence $(f_j)$, with $f_j\in F_j$ 
for all $j$, has a subsequence equivalent to the $\ell_1$-basis. Then 
there is a large weakly null {\rm FDD} refinement $(G_j)$ of $(F_j)$, 
having one of the following mutually exclusive alternatives: 
\vskip1pt
{\rm 1)} $(G_j)$ is uniformly anti-Banach-Saks; that is, 
there is a $\delta>0$ so that 
$$\mathop{\underline{\lim}}_{n\to\infty} {1\over n} 
\Big\| \sum_{j=1}^n g_j \Big\| \ge \delta 
\mathop{\underline{\lim}}_{n\to\infty} \|g_n\|$$ 
for all strictly increasing sequences $(n_j)$ in $\nat$ and all 
sequences $(g_j) \in \prod_{j=1}^\infty Ba\, G_{n_j}$. 
\vskip1pt 
{\rm 2)} $(G_j)$ is uniformly Banach-Saks; that is, there is a sequence 
$(\varep_j)$ of positive numbers tending to zero so that 
$${1\over n} \Big\| \sum_{j=1}^n g_j\Big\| \le\varep_n\ \hbox{ for all }
\ n\ ,$$ 
all strictly increasing sequences $(n_j)$ in $\nat$, and all sequences 
$(g_j) \in \prod_{j=1}^\infty Ba\, G_{n_j}$.
Moreover if the $F_n$'s have bases $\bb_n$ with uniformly bounded 
basis constants, $(G_n)$ may be chosen to be a block refinement of 
$(F_n)$ with respect to $(\bb_n)$.\smallskip

\beginsection{\S2. Proofs of Theorems 1 and 6.}

{\it Proof of Theorem 1\/}. 
Without loss of generality we can assume that $X= C(K)$, the space of all 
real or complex valued continuous functions on a compact metric space $K$. 
For $f\in C(K)$ we let $f^+ = \max (f,0)$ in the real case; in the 
complex case we put $f^+ = \min ((\Re f)^+, (\Im f)^+)$. For $A\subset K$ 
we let $\|\cdot\|_A$ be the seminorm on $C(K)$ defined by  $\|f\|_A 
= \sup_{\xi\in A} |f(\xi)|$. Let $(F_n)$ be a large sequence of finite 
dimensional subspaces of $C(K)$. Since $(F_n)$ {\it has\/} a large 
FDD-refinement, we assume without loss of generality that $(F_n)$ 
is already an FDD. 

We consider the following two cases. 
\medskip

\noindent {\caps Case 1:}

\item{(1)} For all nonempty closed sets $\tilde K\subset K$, all $\varep >0$ 
and all large refinements $(H_n)$ of $(F_n)$ there is a relatively open 
set $U\subset \tilde K$, $U\ne\emptyset$, and a large refinement $(\tilde H_n)$ 
of $(H_n)$ so that 
$$\sup\|f\|_U < \varep\ ,\ \hbox{ for  }\ f\in \bigcup_{n\in\nat} 
S_{\tilde H_n}\ .$$

\noindent {\caps Case 2:} 

\item{(2)} There are a nonempty closed set $K_0 \subset K$, 
$\varep_0>0$ and a large refinement $(H_n)$ of $(F_n)$ so that for all 
nonempty and relatively open sets  $U\subset K_0$ and all further large 
refinements $(\tilde H_n)$ of $(H_n)$, 
$$\liminf_{n\to\infty}\ \sup_{h\in S_{\tilde H_n}} \|h\|_U > \varep_0\ .$$ 
\medskip

Clearly, cases 1 and 2 are mutually exclusive and the failure of one implies 
the other holds. We will show that assuming case~1, we can find a weakly 
null large refinement $(G_n)$ of $(F_n)$. Assuming case~2, we shall 
produce a uniformly-$\ell_1$ large refinement $(G_n)$ of $(F_n)$. 

Assume that (1) is satisfied and let $\varep>0$ be arbitrary. Let $K^{(0)}
=K$ and $(H_n^{(0)}) = (F_n)$. We will choose by transfinite 
induction for each $\alpha <\omega_1$ (where $\omega_1$ is the first 
uncountable ordinal), a closed subset $K^{(\alpha)}$ of $K$ and a large 
refinement $(H^{(\alpha)}_n)$ of $(F_n)$, so that 

\item{(3)} $K^{(\beta)} \subseteq K^{(\alpha)}$ and, if $K^{(\alpha)} 
\ne\emptyset$, 
then $K^{(\beta)} \subsetneqq K^{(\alpha)}$, whenever $\alpha<\beta$.

\item{(4)} Except for perhaps finitely many elements, $(H^{(\beta)}_n)$ is a 
refinement of $(H^{(\alpha)}_n)$ whenever $\alpha<\beta$. 

\item{(5)} For all $\xi \in K\setminus K^{(\alpha)}$,  
$$\limsup_{n\to\infty} \ \sup_{f\in S_{H^{(\alpha)}_n}} |f(\xi)| \le
\varep$$

\noindent 
Assume that for some $\alpha <\omega_1$, $(K^{(\gamma)})_{\gamma <\alpha}$ 
and $(H_n^{(\gamma)})_{\gamma <\alpha}$ have been chosen. 
If $\alpha = \gamma+1$ and $K^{(\gamma)} = \emptyset$ set $K^{(\alpha)} = 
\emptyset$ and $(H_n^{(\alpha)}) = (H_n^{(\gamma)})$. 
If $\alpha = \gamma +1$ and $K^{(\gamma)} \ne\emptyset$, by (1) there 
exists a large refinement $(H_n^{(\alpha)})$ of $(H_n^{(\gamma)})$ and a 
relatively open set $U\subset K^{(\gamma)}$, $U\ne\emptyset$, so that 
$$\|f\|_U < \varep\ \hbox{ for all }\ f\in \bigcup_{n\in\nat} 
S_{H_n^{(\gamma)}}\ .$$ 
Set $K^{(\alpha)} = K^{(\gamma)} \setminus U$. 

If $\alpha = \lim_{n\to\infty} \gamma_n$ for some strictly increasing 
sequence $(\gamma_n)$, set $K_\alpha = \bigcap_{n\in\nat} K_{\gamma_n}$ 
and let $(H_n^{(\alpha)})$ be a ``diagonal sequence'' of $(H_n^{(\gamma_m)})_
{n,m\in\nat}$, chosen such that for each $m$, except for perhaps finitely 
many terms, $(H_n^{(\alpha)})$ is a large refinement of $(H_n^{(\gamma_m)})$. 

Since $K$ is compact and metric, (thus $K$ satisfies the Lindel\"off 
condition) we conclude that for some $\alpha <\omega_1$, $K^{(\beta)} = 
K^{(\alpha)}$ for $\alpha \le \beta  <\omega_1$. By (3) it follows that 
$K^{(\alpha)} =\emptyset$ and from (5) it follows that for all $\xi\in K$, 
$$\limsup_{n\to\infty}\ \sup_{f\in S_{H_n^{(\alpha)}}} 
|f(\xi)| \le\varep\ .$$ 
We let $(H_n^{(\varep)}) := (H_n^{(\alpha)})$. 
Repeating this argument for a sequence $(\varep_m)\subset \real_+$ with 
$\varep_m\downarrow 0$ one obtains for each $m\in\nat$, a large refinement  
$(H_n^{(\varep_m)})_{n\in\nat}$, of 
$(H_n^{(\varep_{m-1})})$, satisfying 
for all $\xi \in K$, 
$$\limsup_{n\to\infty}\ \sup_{f\in S_{H_n^{(\varep_m)}}} |f(\xi)| 
\le \varep_m\ .$$ 
If we let $(G_n)$ be a diagonal sequence of $(H_n^{(\varep_m)})_{n,m\in\nat}$, 
still satisfying $\lim_{n\to\infty} \dim (G_n) =\infty$, we deduce that 
for all $\xi\in K$, 
$$\lim_{n\to\infty}\ \sup_{f\in S_{G_n}} |f(\xi)| =0\ .$$
Thus $(G_n)$ is a weakly null large refinement of $(F_n)$. 

We now assume that (2) is satisfied and let $K_0\subset K$, $\varep_0>0$ 
and $(H_n)$ be as in (2). Let $\varep_1 = \varep_0$ in the 
real case and $\varep_1 = \varep_0/\sqrt2$ in the complex case. Let $D$ be 
a countable dense subset of $K_0$. By passing to a large refinement of 
$(H_n)$ we can assume that 
$$\lim_{n\to\infty} \ \sup_{f\in S_{H_n}} |f(\xi)| = 0\ \hbox{ for all }\ 
\xi\in D\ .\leqno(6)$$
Indeed, let $d_1,d_2,\ldots$ be an enumeration of $D$ and let 
$m_1<m_2<\cdots$ be such that $\dim H_{m_n} \ge 2n$; then set $H'_n = 
\{ x\in H_{m_n} : x(d_i) = 0$ for $1\le i\le n\}$. Now 
$\dim H'_n \ge n$ for all $n$, so $(H'_n)$ is the desired large 
refinement. 
Let $\varep_1/34 > \delta >0$. By induction we will choose 
an increasing sequence of integers $(k_n)$ and for each $n$, 
a subspace $G_n$ of $H_{k_n}$ and a finite set $\Pi_n$ consisting of 
nonempty relatively open subsets of $K_0$ so that the following conditions are 
satisfied: 
\smallskip
\item{(7)} $\dim (G_n) \ge n$, 
\smallskip
\noindent and
\smallskip
\item{(8)} For every $g\in S_{G_n}$, and every $U \in \Pi_{n-1}$ (let 
$\Pi_0 = \{ K_0\}$) there are $U_1,U_2 \in \Pi_n$, $U_1\cup U_2\subseteq U$, 
so that 
$$g^+|_{U_1} \ge \varep_1 -\delta \quad \hbox{and}\quad 
\|g\|_{U_2} \le \delta\ .$$

\noindent Once  we have chosen $(G_n)$ in this way we conclude that $(G_n)$ 
must be  uniformly-$\ell_1$. To see this, fix $(f_n)$ with $f_n\in S_{G_n}$ 
for all $n\in \nat$. For each $n$, let $A_n= \{ k\in K: f_n (k) 
>\varep_1-\delta\}$ 
and $B_n = \{k\in K : |f_n(k)| <\delta\}$. Evidently $A_n\cap B_n = 
\emptyset$ for all $n$. We shall show that $(A_n,B_n)$ is an independent 
sequence of pairs, in the terminology of [R1]. Once this is done, a 
refinement of the argument in [R1] yields that $(f_n)$ is 
$\frac{16}{\varep_1}$-equivalent to the $\ell_1$-basis. 

Indeed, we first can inductively choose 
sets $(U_i^{(n)} : i=1,2,\ldots,2^n)\subset \Pi^{(n)}$ so that 
$$f_n^+\big|_{\bigcup\limits_{i=1}^{2^{n-1}} U_{2i-1}^{(n)}} > \varep_1
-\delta\quad\hbox{and}\quad 
\| f_n\|_{\bigcup\limits_{i=1}^{2^{n-1}} U_{2i}^{(n)}} < \delta$$ 
and so that $U_{2j}^{(n)} \cup U_{2j-1}^{(n)} \subset U_j^{(n-1)}$ for 
$n\in\nat$ and $j=1,2,\ldots,2^{n-1}$. 
Now fix $N$, $I$ and $J$ non-empty disjoint subsets of $\{1,\ldots,
N\}$, say with $I\cup J = \{1,\ldots,N\}$. We see that $\bigcap_{n\in I} 
A_n \cap \bigcap_{n\in J} B_n$ is non-empty by defining the following 
sequence of sets $C_0,C_1,\ldots,C_N$: Let $U_1^0 = K_0 = C_0$, 
$1\le n\le N$, and suppose $C_{n-1}$ is chosen with $C_{n-1} = U_j^{(n-1)}$ 
for some $1\le j\le 2^{n-1}$. If $n\in I$, set $C_n = U_{2j-1}^{(n)}$, 
otherwise set $C_n = U_{2j}^{(n)}$. Then the $C_n$'s satisfy that 
$\bigcap_{n=1}^N C_n \ne \emptyset$ and for all $n$, $C_n\subset  A_n$ 
if $n\in I$, $C_n\subset B_n$ if $n\in J$. 

Now let $\sum_{j=1}^N |a_j| =1$ with $a_j = b_j + ic_j$ for $j\le N$. 
By multiplying by $-1$, $i$ or $-i$ if necessary we may assume that 
$\sum_{j=1}^N b_j^+ \ge 1/4$. Let $I = \{ j\le N: b_j \ge 0$ and 
$c_j \ge 0\}$ and $J= \{j\le N: b_j \ge 0$ and $c_j<0\}$. Thus either 
$$\sum_{j\in I} (b_j + c_j) \ge {1\over8}\qquad\hbox{or}\qquad
\sum_{j\in J} (b_j -c_j) \ge {1\over8}\ .$$ 
Suppose the first sum exceeds $1/8$. 
Now by the independence of $(A_n,B_n)$, 
choose $k\in K$ such that 
$f_j^+ (k) >\varep_1 -\delta$ for $j\in I$ and 
$|f_j(k)| <\delta$ for $j\notin I$. Let $f_j(k) = B_j +iC_j$. Then 
$$\eqalign{ 
\Big| \sum_{j=1}^N a_j f_j(k)\Big| 
& \ge \Big| \Im \biggl( \sum_{j=1}^N a_j f_j(k)\biggr)\Big|\cr 
\myskip 
& = \Big| \sum_{j=1}^n (b_j C_j + B_j c_j)\Big|\cr 
\myskip 
& \ge \sum_{j\in I} (b_j C_j + B_j c_j) 
- \sum_{j\notin I} |b_j C_j + B_j c_j|\cr 
\myskip 
&\ge {(\varep_1-\delta)\over 8} - 2\delta > {\varep_1\over16}\ .\cr}$$
A similar estimate ensues if the second sum exceeds $1/8$.  
Thus $(f_n)$ is indeed $16\over\varep_1$-equivalent to the $\ell^1$-basis. 

Assume that for some $n\ge 1$, $\Pi_{n-1}$ and $k_{n-1}$ (let $k_0=0$) 
are chosen. 
Now consider the finite family of Lipschitz functions defined on 
$C(K)$ by $f\mapsto \|f^+\|_U$, $U\in \Pi_{n-1}$. Since $(H_n)$ is large, 
we may use the first stabilization principle in order to 
pass to a large refinement $(\tilde H_i)$ of $(H_i)_{i>k_{n-1}}$ so that 
for some family ($a_i^{(U)}: U\in\Pi_{n-1}$, $i\in\nat$) in $\real^+$ we have 
$$a^{(U)}_i - {\delta\over 4} < \|f^+ \|_U < a^{(U)}_i + {\delta \over 4}$$ 
whenever $U\in \Pi_{n-1}$, $i\in \nat$ and $f\in S_{\tilde H_i}$. 
>From (2) we deduce that there exists $i_0\in\nat$ so that for all $i\ge i_0$ 
and $U\in \Pi_{n-1}$ we have $a_i^{(U)} \ge \varep_1 - {\delta\over4}$. 
Indeed, in the real case we only have to observe that if $\|f\|_U \ge 
\varep_0$ then $\| f^+\|_U \ge \varep_0$ or $\| (-f)^+\|_U \ge\varep_0$; 
in the complex we find for any $f\in C(K)$ for which $\|f\|_U > \varep_0$,  
a point $\xi\in U$ with $|f(\xi)|>\varep_0$ and then a complex number $a$, with 
$|a| =1$, so that $\Re (a\cdot f(\xi)) = \Im (a\cdot f(\xi)) = (a\cdot f
(\xi))^+$. Thus $\|( a\cdot f)^+\|_U\ge {1\over\sqrt2} \|f\|_U  > \varep_1$. 
We deduce that 
$$\| f^+\|_U > \varep_1 - {\delta\over2}\leqno(9)$$
for all $U\in\Pi_{n-1}$, $i\ge i_0$ and $f\in S_{\tilde H_i}$. 

Now using (6), 
we pick, for each $U\in \Pi_{n-1}$, an element $\xi_U\in U\cap D$ and find an 
$i_1 \ge i_0$ so that $\dim (\tilde H_{i_1})\ge n$ and so that 
$$\sup_{f\in S_{\tilde H_{i_1}}} |f(\xi_U)| < {\delta \over 2}\ . 
\leqno(10)$$ 
Let $(f_s)_{s=1}^\ell$ be a finite  ${\delta\over2}$-net for 
$S_{\tilde H_{i_1}}$. We find by (9) and (10) for each $U\in \Pi_{n-1}$, 
non-empty open subsets $V_0^{(U)},V_1^{(U)},\ldots,V_\ell^{(U)}$ so that 
$f_s^+|_{V_s^{(U)}} > \varep_1 -{\delta\over2}$ 
and $\|f_s\|_{V_0^{(U)}} < {\delta\over2}$, 
for $s=1,2,\ldots,\ell$.  
This implies that for all 
$f\in S_{\tilde H_{i_1}}$ we have $\|f\|_{V_0^{(U)}} <\delta$, and for some 
$1\le s\le \ell$ (namely the $s$ for which $\|f-f_s\| <{\delta\over2}$) 
we have $f^+ |_{V_s^{(U)}} > \varep_1 -\delta$. 
Set 
$ \Pi_n = \left\{ V_0^{(U)} : U\in \Pi_{n-1}\right\} 
\cup \left\{ V_s^{(U)} : 1\le s\le \ell\ ,\ U\in \Pi_{n-1}\right\}$,
$G_n=\tilde H_{i_1}$, and choose 
$k_n > k_{n-1}$ so that $\tilde H_{i_1} \subset 
H_{k_n}$. This completes the induction and thus the proof of the first 
version of Theorem~1. 

The ``block-version'' of Theorem 1 is proved in exactly the same way using 
the second stabilization principle instead of the first. One need only 
note that block refinements could be taken wherever we took simple 
refinements.\qed 

\noindent {\it Proof of Theorem 6\/}. 
As in the proof of Theorem 1 we will only show the first version of 
Theorem~6. The ``block-version'' is left to the reader. We shall assume 
that $X$ is a Banach space over $\real$. The complex case does not 
provide any further difficulties. 

Let $f:X\to \real$ be Lipschitz and let $\varep_n\downarrow 0$. 
We accomplish the proof by induction, insuring the conditions in the 
Theorem for a fixed $k\ge 2$. Precisely, we shall choose for each $k$,  
a large sequence $(G_n^{(k)})$ of finite dimensional subspaces so that  
$(G_n^{(k+1)})$ is a refinement of $(G_n^{(k)})$ ($(G_n^{(1)})\equiv (F_n)$), 
and a function $C^{(k)} :Ba (\ell_\infty^k) \to\real$, so that 
$$f(\alpha_1x_1 +\cdots + \alpha_k x_k) \buildrel {\varep_{n_1}}\over = 
C^{(k)} (\alpha_1,\ldots,\alpha_k)$$ 
whenever $(\alpha_1,\ldots,\alpha_k) \in Ba (\ell_\infty^k)$ and 
$x_k \in S(G_{n_i}^{(k)})$, for all $1\le n_1<n_2< \cdots < n_k$. 

Once this is done, then by diagonalization 
we finally find a large refinement $(G_n)$ of $(F_n)$ and functions 
$C^{(k)} :Ba (\ell_\infty^k)\to\real$, $k\in\nat$ so that for all $k\in\nat$ 
and all $k\le n_1 < n_2 < \cdots < n_k$ we have 
$$f(\alpha_1 x_1 +\alpha_2 x_2 +\cdots + \alpha_k x_k) 
\buildrel {\varep_{n_1}}\over = 
C^{(k)} (\alpha_1,\ldots,\alpha_k)$$ 
whenever $x_i \in S_{G_{n_i}}$, for $i=1,2,\ldots,k$. 

Clearly we have that 
$$C^{(k)}(\alpha_1,\alpha_2,\ldots,\alpha_k) = C^{(k+s)} 
(\alpha_1,\alpha_2,\ldots,\alpha_k,0,0,\ldots,0)$$ 
for $k,s\in \nat$ and $(\alpha_1,\alpha_2,\ldots,\alpha_k) 
\in Ba(\ell_\infty^k)$, and, thus, if we put 
$$\tilde f(\alpha_1,\ldots,\alpha_k,0,0,\ldots) = 
C^{(k)} (\alpha_1,\alpha_2,\ldots,\alpha_k)\ ,$$ 
for $k\in \nat$ and $(\alpha_i)_{i=1}^k \in Ba(\ell_\infty^k)$, 
$\tilde f$ has the required properties.

We now indicate in detail how to carry this out for $k=2$. 
First note the following 
\medskip

\noindent {\it Fact\/}. 
Let $g: S_X\to \real$ be Lipschitz and let $(L_n)$ be any large sequence 
of finite dimensional subspaces of $X$. Let $\delta_n\downarrow 0$. 
There exist a large refinement $(\tilde L_n)$ of $(L_n)$ and 
$C\in\real$ such that for all $n$ and $y\in S_{\tilde L_n}$, 
$$g(y) \buildrel {\delta_n}\over = C\ .$$ 

This follows easily from the first stabilization theorem. One first obtains 
a large refinement $(\tilde{\tilde L}_n)$ of $(L_n)$ and $(C_n)\subseteq 
\real$ such that $g(y) \buildrel {\delta_n/2}\over = C_n$ for 
$y\in S_{\tilde{\tilde L}_n}$. $(C_n)$ is bounded so for some subsequence 
$(C_{k_n})$ and $C\in \real$, $|C_{k_n} - C| < \delta_n/2$ for all $n$. 
Let $\tilde L_n = \tilde{\tilde L}_{k_n}$. 

Let $H_1\equiv F_1$. Choose finite sets $D_1\subseteq D_2\subseteq \cdots 
\subseteq Ba (\ell_\infty^2)$ and $\D_1 \subseteq \D_2\subseteq \cdots 
\subseteq S_{H_1}$ so that for all $n$, $D_n$ is an $\varep_n$-net for 
$Ba(\ell_\infty^2)$ and $\D_n$ is an $\varep_n$-net for $S_{H_1}$. 
For $x\in S_{H_1}$ and $(\alpha,\beta) \in Ba (\ell_\infty^2)$, 
$y\mapsto f(\alpha x+\beta y)$ is a Lipschitz function on $X$. Thus by 
iterating the fact above a finite number of times we obtain a large 
refinement $(F_n^{(1,1)})_{n=1}^\infty$ of $(F_n)$ and 
$(C(\alpha,\beta,x))_{(\alpha,\beta,x)\in D_1\times \D_1}\subseteq \real$ 
such that for all $(\alpha,\beta) \in D_1$, $x\in \D_1$ and 
$y\in F_n^{(1,1)}$, 
$$f(\alpha x+\beta y) \buildrel {\varep_n}\over = C(\alpha,\beta,x)\ .$$ 
Repeating this argument inductively we obtain for all $k\in \nat$, a large 
refinement $(F_n^{(1,k)})_{n=1}^\infty$ of $(F_n^{(1,k-1)})_{n=1}^\infty$ 
and $(C(\alpha,\beta,x))_{(\alpha,\beta,x)\in D_k\times \D_k}$ such that 
$$f(\alpha x + \beta y) \buildrel {\varep_n}\over = C(\alpha,\beta,x)$$ 
if $(\alpha,\beta) \in D_k$, $x\in \D_k$ and $y\in F_n^{(1,k)}$. 
By diagonalization we obtain a large refinement $(F_n^{(1)})_{n=1}^\infty$ 
of $(F_n)$ with the property 
\smallskip
\iitem{i)} For $k\in\nat$, $(\alpha,\beta,x) \in D_k\times \D_k$, 
$n\ge k$, and $y\in F_n^{(1)}$, 
$$f(\alpha x+\beta y) \buildrel {\varep_n}\over = C(\alpha,\beta,x)\ .$$ 
\smallskip

\noindent Suppose that the Lipschitz constant of $f$ is $K\ge 1$, 
i.e., $|f(x)-f(y)| \le K\|x-y\|$. Then for $(\alpha,\beta),(\alpha',\beta') 
\in Ba(\ell_\infty^2)$, $x,x'\in S_{H_1}$, and $\|y\|=1$, we have 
$$\eqalign{
|f(\alpha x+\beta y) - f(\alpha'x'+ \beta' y)| 
&\le K\|(\alpha x - \alpha' x) + (\alpha' x-\alpha'x') 
+ (\beta -\beta')y\| \cr 
&\le K\bigl( |\alpha-\alpha'| +  |\beta-\beta'| + \|x-x'\|\bigr)\ .\cr}$$ 
>From this and i) we obtain 
\smallskip
\iitem{ii)}\quad  $|C(\alpha,\beta,x) - C(\alpha',\beta',x')| 
\le K\bigl[ |\alpha-\alpha' | + |\beta-\beta'| + \|x-x'\|\bigr]$ 
\smallskip

\noindent whenever $(\alpha,\beta),(\alpha',\beta') \in \cup D_n$ and 
$x,x'\in \cup \D_n$. Thus we can uniquely extend $C(\alpha,\beta,x)$ 
to a function $C^{[1]}: Ba (\ell_\infty^2 ) \times S_{H_1}\to \real$ 
which satisfies ii) for all $(\alpha,\beta),(\alpha',\beta') \in 
Ba(\ell_\infty^2)$ and $x,x'\in S_{H_1}$. Furthermore, by replacing 
$(F_n^{(1)})$ by an appropriate subsequence, we may assume that 
\smallskip
\iitem{iii)} For $n\in \nat$, $(\alpha,\beta,x)\in Ba(\ell_\infty^2)
\times S_{H_1}$ and $y\in S_{F_n^{(1)}}$, 
$$f(\alpha x+\beta y) \buildrel {\varep_n}\over = 
C^{[1]} (\alpha,\beta,x)\ .$$ 

Set $H_2 = F_{n_2}^{(1)}$ where $n_2$ is chosen so that $\dim H_2 >\dim H_1$. 
Proceeding as above we obtain a function 
$$C^{[2]} : Ba (\ell_\infty^2) \times S_{H_2} \to \real$$ 
and a large refinement $(F_n^{(2)})$ of $(F_n^{(1)})$ so that 
\smallskip
\iitem{iv)} \quad $|C^{[2]} (\alpha,\beta,x) - C^{[2]}(\alpha',\beta',x')| 
\le K\bigl[ |\alpha-\alpha'| + |\beta-\beta'| + \|x-x'\|\bigr]$ 
\smallskip

\noindent for all $(\alpha,\beta),(\alpha',\beta') \in Ba (\ell_\infty^2)$ 
and $x,x'\in S_{H_2}$ and 
\smallskip
\iitem{v)} \quad $f(\alpha x+\beta y) \buildrel {\varep_n}\over = 
C^{[2]} (\alpha,\beta,x)$ 
\smallskip

\noindent for all $x\in S_{H_2}$, $(\alpha,\beta) \in Ba (\ell_\infty^2)$ 
and $y\in S_{F_n^{(2)}}$. 

We continue in this manner obtaining a large refinement $(H_n)$ of 
$(F_n)$ and, for $k\in \nat$, functions $C^{[k]}: Ba (\ell_\infty^2) 
\times S_{H_k}\to\real$ satisfying 
\smallskip
\iitem{vi)} \quad $|C^{[k]} (\alpha,\beta,x) - C^{[k]} (\alpha',\beta',x')| 
\le K\bigl[ |\alpha-\alpha'| + |\beta-\beta'| + \|x-x'\|\bigr]$ 
\smallskip

\noindent for $(\alpha,\beta),(\alpha',\beta')\in Ba (\ell_\infty^2)$ 
and $x,x'\in S_{H_k}$ and 
\smallskip
\iitem{vii)} \quad $f(\alpha x+\beta y) \buildrel {\varep_n}\over = 
C^{[k]} (\alpha,\beta,x)$ 
\smallskip

\noindent for all $(\alpha,\beta) \in Ba (\ell_\infty^2)$, $x\in S_{H_k}$ 
and $y\in S_{H_n}$ with $n>k$. 
\medskip

\noindent (Actually it might be necessary to pass to a subsequence of 
$(H_n)$ to obtain the precise estimate vii).) 
\smallskip

We now apply the first stabilization result to finite sets of functions 
$C^{[k]} (\alpha,\beta,\cdot)$. Let $n\in \nat$ and $0\le \bar\varep \le 
\min \{|\alpha-\alpha'| + |\beta -\beta'| : (\alpha,\beta) \ne (\alpha',
\beta') \in D_n\}$. Consider for a fixed $k$ the Lipschitz functions 
$C^{[k]} (\alpha,\beta,\cdot) : S_{H_k} \to \real$ for $(\alpha,\beta) 
\in D_n$. If $\dim H_k$ is sufficiently large there exists $\tilde H_k 
\subseteq H_k$, $\dim \tilde H_k \ge n$ and $(C^k(\alpha,\beta))_{(\alpha,
\beta)\in D_n} \subseteq \real$  
so that $C^{[k]} (\alpha,\beta,x) 
\buildrel {\bar \varep}\over = C^k (\alpha,\beta)$ for all 
$x\in S_{\tilde H_k}$ and $(\alpha,\beta) \in D_n$. Thus this plus vi) yields 
$$\eqalign{|C^k (\alpha,\beta) - C^k (\alpha',\beta') | 
&\le K\bigl[ |\alpha-\alpha' | + |\beta-\beta'|\bigr] +2\bar\varep\cr 
& \le 3K\bigl[ |\alpha-\alpha'| + |\beta-\beta'|\bigr]\cr}$$ 
for $(\alpha,\beta),(\alpha',\beta') \in D_j$. 
The last inequality holds by the choice of $\bar\varep$ and the fact 
that $K\ge1$. 

We inductively use this argument for the parameters $(n,\bar\varep_n)$ 
where $\bar\varep_n\downarrow 0$ rapidly chosen depending on $(D_n)$ and 
$(\varep_n)$. We obtain a large refinement $(I_n)$ of $(H_n)$ with 
${\dim I_n\ge n}$ and functions $C^n :D_n \to\real$ satisfying 
\smallskip
\iitem{viii)} \quad $|C^n (\alpha,\beta) - C^n (\alpha',\beta')| 
\le 3K\bigl[ |\alpha-\alpha'| + |\beta-\beta'|\bigr]$ for $(\alpha,\beta), 
(\alpha',\beta') \in D_n$ 
\smallskip
\noindent and 
\smallskip
\iitem{ix)} For all $x\in S_{I_n}$ and $(\alpha,\beta) \in D_n$, 
$$f(\alpha x+\beta y) \buildrel {\varep_n}\over = 
C^n(\alpha,\beta)\ \hbox{ whenever }\ y\in S_{I_q}\ ,\qquad q> n\ .$$  

\noindent 
For  $n\in \nat$, the function $(C^k |_{D_n})_{k\ge n}$ are uniformly 
Lipschitz. Thus by a compactness argument we can find a Lipschitz function 
$C^{(2)} : \cup D_j\to \real$ and $k_1<k_2<\cdots$ so that for all $n$
and $(\alpha,\beta)\in D_n$, 
$$C^{(2)} (\alpha,\beta) \buildrel {\varep_n}\over = 
C^{k_n} (\alpha,\beta)\ .$$ 
$C^{(2)}$ thus uniquely extends to a continuous function $C^{(2)} :Ba 
(\ell_\infty^2) \to \real$. Letting $(G_n^{(2)})_{n=1}^\infty$ 
be a suitable subsequence of $(I_{k_n})$ we obtain 
\smallskip
\iitem{x)} \quad $f(\alpha x_{n_1} + \beta x_{n_2}) 
\buildrel {\varep_{n_1}}\over = C^{(2)} (\alpha,\beta)$ for all 
$(\alpha,\beta) \in Ba (\ell_\infty^2)$, 
\iitem{} \quad $n_1<n_2$, $x_{n_1}\in 
S_{G_{n_1}^{(2)}}$ and $x_{n_2} \in S_{G_{n_2}^{(2)}}$
\smallskip

\noindent which was what was needed to be proved in the case $k=2$.\qed

\beginsection{\S3. A Sketch of the Proof of the Second Stabilization 
Principle}

The reader unfamiliar with Lemberg's proof might first wish to read that 
argument (see [MS, Ch.12]). 
In order to shorten the proof we will not only use Lemberg's proof of 
Krivine's theorem but also the quantitative version of this theorem. 

\proclaim Theorem 9. {\rm (see [R3])}\quad 
For every $C>1$, $\varep >0$ and $k\in\nat$, there is an $n= n(C,\varep,k) 
\in\nat$ so that: If $F$ is a Banach space of dimension $n$ and if 
$(f_i)_{i=1}^n$ is a basis of $F$ having basis constant not exceeding $C$,  
then there exists a block basis $(g_i)_{i=1}^k$ of $(f_i)_{i=1}^n$ and a 
$p\in [1,\infty]$ so that $(g_i)_{i=1}^k$ is $(1+\varep)$-equivalent to 
the unit basis of $\ell_p^k$. 

In view of Theorem 9 we only have to prove the second stabilization 
principle for finite dimensional $\ell_p$-spaces. Using a compactness 
argument, similar to the argument of [R3] 
by which Theorem~9 was deduced from the 
finite dimensional version of Krivine's theorem, we only have to 
show the following claim. 

\proclaim Claim 1. 
Let $X= \ell_p$, $1\le p<\infty$, or $X= c_0$, and let $f:X\to\real$ be a 
Lipschitz function. For each $\varep >0$ and $k\in\nat$ there exists a 
block basis $(y_i)_{i=1}^k$ of $(e_i)$ (the unit vector basis of $X$) so that 
$\osc (f|_{S_{[y_i]_{i=1}^k}}) < \varep$. 

\noindent {\it Proof of Claim 1.} 
We need some notation. For $x,y\in c_{00}$ we say $x$ and $y$ have the 
{\it same distribution\/}, and write $x\disteq y$, if $x= \sum_{i=1}^k 
\alpha_i e_{n_i}$ and $y= \sum_{i=1}^k \alpha_i e_{m_i}$ for some 
$k\in\nat$, $(\alpha_i)_{i=1}^k \subset \real$, resp. $\complex$, and 
$n_1< n_2< \cdots n_k$ and $m_1< m_2< \cdots m_k$. 
We define for $x,y\in X\cap c_{00}$, 
$$\dis (x,y) = \inf \left\{ \| \bar x-\bar y\| : 	
\bar x\disteq x\ \hbox{ and }\ \bar y \disteq y\right\}\ .$$ 
For $x\in c_{00}$, we let $\supp (x) = \{i\in\nat: x_i\ne0\}$, and 
write $x<y$ for $x,y\in c_{00}$ if $\max (\supp (x)) < \min (\supp(y))$. 

We first reduce claim 1 to the case that $f$ is {\it $1$-unconditional\/} and 
{\it $1$-spreading\/}. By this we mean that $f(\sum \alpha_i e_i) = 
f(\sum |\alpha_i|e_{n_i})$ for all $\sum \alpha_i e_i \in X$ and all 
strictly increasing sequences $(n_i)\subset \nat$. 

In order to reduce claim 1 to the 1-unconditional and 1-spreading case we 
first pass to a sequence $n_i\subset \nat$ for which 
$$f^{(1)} \biggl( \sum_{i=1}^\ell \alpha_i e_i\biggr) 
= \lim_{k_1\to\infty}\ \lim_{k_2\to\infty} \ldots \lim_{k_\ell\to\infty} 
f\biggl( \sum_{i=1}^\ell \alpha_i e_{n_{k_i}}\biggr)$$ 
exists for all $\ell$ and scalars $\alpha_1,\alpha_2,\ldots,\alpha_\ell$. 
It follows that $f^{(1)}$ is 1-spreading on $X$. If $X$ is defined over $\real$ 
we let $\ell_n =2$, for $n\in\nat$, and put $(\xi_i)_{i=1}^{\ell_n} = (1,-1)$. 
If $X$ is defined over $\complex$ we let $\ell_n=n$, for $n\in\nat$, 
and $(\xi_j)_{j=1}^{\ell_n} = e^{(i2\pi j/n)}$. 
If $X=\ell_p$, $1\le p<\infty$, let $(u_n)$ be a sequence in $X$ with 
$u_1 <u_2 < \cdots$ and 
$$u_n \disteq {1\over (n\ell_n)^{1/p}} 
\sum_{s=1}^n\ \sum_{t=1}^{\ell_n} \xi_t e_{(s-1)\ell_n+t}\ ,\ \hbox{ for }\ 
n\in\nat\ .$$ 

If $X= c_0$ we let $(u_n)$ be a sequence in $X$, with $u_1<u_2<\cdots$, and 
$$u_n \disteq \sum_{s=1}^n {s\over n} \cdot 
\sum_{t=1}^{\ell_n} \xi_t e_{(s-1)\ell_n+t} + 
\sum_{s=1}^{n-1} {n-s\over n} \sum_{t=1}^{\ell_n} \xi_t  e_{(n+s-1)\ell_n+t}$$ 
Note that $(u_n)$ is normalized, and that from the fact that $f^{(1)}$ is 
1-spreading it follows that for some sequence $\varep_n\downarrow 0$ and 
some subsequence $(\tilde u_n)$ of $(u_n)$, 
$$f^{(1)} \biggl( \sum_{j=1}^k \alpha_j \tilde u_{j+n}\biggr) 
\buildrel {\varep_n}\over = 
f^{(1)} \biggl( \sum_{j=1}^k \sigma_j \alpha_j \tilde u_{j+n}\biggr)$$ 
whenever $k,n\in\nat$, $|\sigma_j|=1$, for $j=1,2,\ldots,k$, and 
$\|\sum_{j=1}^k \alpha_j e_j\|=1$. 

Pass now to a subsequence $(n_i)\subset \nat$ for which 
$$f^{(2)} \biggl( \sum_{i=1}^\ell \alpha_i e_i\biggr) \equiv 
\lim_{k_1\to\infty}\ \lim_{k_2\to\infty}\ldots \lim_{k_\ell\to\infty} 
f^{(1)} \biggl( \sum_{i=1}^\ell \alpha_i \tilde u_{n_{k_i}}\biggr)$$ 
exists, whenever $\ell\in\nat$ and $(\alpha_i)_{i=1}^\ell \in c_{00}$. 
$f^{(2)}$ is 1-unconditional and 1-spreading and we need only prove that 
claim~1 is true for $f^{(2)}$. 
Thus, in order to finish the proof of claim~1 we need to show the 
following claim~2. 

\proclaim Claim 2. 
For every $\varep >0$ and $k\in\nat$ there is a block basis 
$(x_i)_{i=1}^k$ of $(e_i)$ which is normalized in $X$, having the 
property that the set 
$$B^+ (x_1,\ldots,x_k) = \left\{ \sum_{i=1}^k \alpha_i x_i : 
0\le \alpha_i \le1\ ,\ \Big\| \sum_{i=1}^k \alpha_i x_i\Big\| = 1\right\}$$ 
has diameter less than $\varep$ with respect to $\dis (\cdot,\cdot)$.

\newpage
\noindent {\it Proof of Claim 2\/}. 

\noindent {\caps Case 1:} $X= \ell_p$, $1\le p<\infty$. 

In this case we consider as in [L] the ``rationalized'' version of $\ell_p$, 
i.e., $\ell_p (D) = \{ (x_q)_{q\in D} :\sum_{q\in D} |x_q|^p<\infty\}$ 
where $D= \que \cap (0,1)$. 
Let $(e_q)_{q\in D}$ denote the natural basis of $\ell_p(D)$. 
For $n\in \nat$ define the operator 
$T_n:\ell_p (D) \to \ell_p(D$ by 
$$T_n \biggl( \sum_{q\in D} \alpha_q e_q\biggr) 
= \sum_{j=1}^n\ \sum_{q\in D} \alpha_q e_{(q+j-1)/n}\ .$$
For every $n\in\nat$, $\lambda_n = n^{1/p}$ is an approximate eigenvalue of 
$T_n$ [L] and since $T_n$ and $T_m$ commute for $n,m\in\nat$ one can choose 
for a fixed $m\in\nat$, $m\gg k$, and $\delta >0$ a vector 
$u=\sum_{q\in D} u_q e_q \in Ba(\ell_p(D))$ so that $\supp (u) = \{q\in D: 
u_q\ne0\}$ is finite, and so that $\|T_n(u) - n^{1/p}u\| <\delta$ 
for all $m\le n$. 

Let $x_1<x_2<\cdots < x_m$ be elements of $\ell_p$ ($=\ell_p(\nat)$), 
each having the same distribution as $u$ (i.e., 
$x_k \disteq \sum_{i=1}^s u_{q_i} e_i$ where $q_1<q_2<\cdots < q_s$ and 
$\supp (u) = \{q_1,q_2,\ldots,q_s\}$). We deduce that for any scalars 
$\alpha_1,\alpha_2,\ldots,\alpha_k$ with $\|\sum_{i=1}^k \alpha_i e_i\|=1$ 
we have 
$$\eqalign{ 
\dis \biggl( x_1\ ,\ \sum_{i=1}^k \alpha_i x_i\biggr) 
&\buildrel {\delta^1}\over = \dis\biggl( x_1\ ,\  \sum_{i=1}^k 
\Bigl( {m_i\over m}\Bigr)^{1/p} x_i\biggr) \cr 
\myskip 
&\buildrel \delta\over \le \dis \biggl( {1\over m^{1/p}} 
\sum_{i=1}^m x_i\ ,\ \sum_{i=1}^k \Bigl( {m_i\over m}\Bigr)^{1/p} x_i\biggr)\cr 
\myskip 
&\le {1\over m^{1/p}} \sum_{i=1}^k \dis  
\bigg( \sum_{j=1}^{m_i} x_j\ ,\ (m_i)^{1/p} x_i\biggr) 
\le k\cdot \delta\ ,\cr}$$ 
where $m_1,\ldots,m_k\in\nat$ with $\sum_{i=1}^k m_i =m$ are chosen so that 
$\sum |\alpha_i - ({m_i\over m})^{1/p}|$ is minimal, and 
where $\delta_1$ depends on $m$ and decreases to zero for $m\to\infty$. 
Thus, choosing $m$ big enough and $\delta$ small enough we deduce 
claim~2 in the case that $X=\ell_p$.
\medskip

\noindent {\caps Case 2:} $X=c_0$ 

In this case Lemberg's argument does not work, but we are able to explicitly 
write down the desired vectors $x_1,x_2,\ldots,x_k$. 

For $0<r<1$ we will define a sequence of vectors $(y^{(n)}: n \in\nat_0)$ 
in $Ba(c_0) \cap c_{00}$. We put $y^{(0)} = e_1$ and assuming 
$y^{(n)} = \sum_{i=1}^{\ell_n} y_i^{(n)} e_i$ is chosen we put 
$$y^{(n+1)} = \sum_{i=1}^{\ell_n} \left( r^{n+1} e_{3(i-1)+1} 
+ y_i^{(n)} e_{3(i-1)+2} + r^{n+1} e_{3(i-1)+3}\right)$$ 
(thus  $y^{(1)} = (r,1,r,0,\ldots)$, 
$y^{(2)} = (r^2,r,r^2,r^2,1,r^2,r^2,r,r^2,0,\ldots)$, etc.). 
Choosing $r<1$ big enough and $n\in\nat$ big enough and letting 
$x_1 < x_2<\cdots < x_k$ all have the same distribution as $y^{(n)}$ one 
also deduces claim~2.\qed 

\noindent {\it Remark.} 
For the case $X=c_0$,  T.~Gowers [G] independently obtained a deeper 
version of claim~1. He showed that  for every Lipschitz function 
$f:c_0\to\real$ and every $\varep>0$ there is an infinite dimensional 
subspace $Y$ of $c_0$ so that $\osc (f|_{S_Y}) <\varep$. 
This is false for $X=\ell_p$ $(1\le p<\infty)$ [OS]. 
\bigskip 

\baselineskip=12pt
\frenchspacing 
\leftline{\bf References}

\ref {[BL]} B. Beauzamy and J.-T. Laprest\'e, 
{\it Mod\`eles \'etal\'es des espaces de Banach}, 
Travaux en Cours, Herman, Paris, 1984. 

\ref {[CGJ]} B.J. Cole and T. Gamelin and W.B. Johnson, 
{\it Analytic Disks in Fibers over the unit ball of a Banach space}, 
Mich. J. Math. (to appear). 

\ref {[G1]} T. Gowers, 
{\it Lipschitz functions on classical spaces}, 
preprint. 

\ref {[G2]} T. Gowers, 
{\it A space not containing $c_0$, $\ell_1$ or a reflexive subspace}, 
preprint. 

\ref {[K]} J.L. Krivine, 
{\it Sous-espaces de dimension finie des espaces de Banach reticul\'es}, 
Ann. of Math. {\bf104} (1976), 1--29. 

\ref {[L]} H. Lemberg, 
{\it Nouvelle d\'emonstration d'un th\'eor\`eme de J.L.~Krivine sur la 
finie repr\'esen-tation de $\ell_p$ dans un espace de Banach}, 
Israel J. Math. {\bf39} (1981), 341--348. 

\ref {[LT]} J. Lindenstrauss and L. Tzafriri, 
``Classical Banach Spaces I,''
Springer-Verlag, New York, 1975. 

\ref {[MS]} V.D. Milman and G. Schechtman, 
``Asymptotic Theory of Finite  Dimensional Normed Spaces,'' 
LNM 1200, Springer-Verlag, New York, 1986. 

\ref {[O]} E. Odell, 
{\it Applications of Ramsey theorems to Banach space theory}, 
in ``Notes in Banach Spaces'' (H.E.~Lacey, ed.), 
University Press, Austin and London, 379--404, 1980.

\ref {[OS]} E. Odell and Th. Schlumprecht, 
{\it The Distortion Problem}, 
preprint. 

\ref {[R1]} H. Rosenthal, 
{\it A characterization of Banach spaces containing $\ell_1$}, 
Proc. Nat. Acad. Sci. USA {\bf71} (1974), 2411--2413. 

\ref {[R2]} H. Rosenthal, 
{\it Weakly independent sequences and the Banach-Saks property}, 
Abstract for the 1975 Durham Symposium, 
Bull. London Math. Soc. {\bf8}  (1976), 22--24. 

\ref {[R3]} H. Rosenthal, 
{\it On a theorem of J.L.~Krivine concerning block finite-representability 
of $\ell_p$ in general Banach spaces}, 
J.~Funct. Anal. {\bf28} (1978), 197--225. 

\ref {[R4]} H. Rosenthal, 
{\it Some remarks concerning unconditional basic sequences}, 
Longhorn Notes (1982-83), The University of Texas at Austin, 15--48. 

\ref {[R5]} H. Rosenthal, 
{\it The unconditional basis sequence problem}, 
Contemp. Math. {\bf52} (1986), 70--88. 

\ref {[R6]} H. Rosenthal, 
{\it Some aspects of the 
subspace structure of infinite dimensional Banach spaces}, 
in ``Approximation Theory and Functional Analysis'' (C.K.~Chui, ed.), 
Academic Press, San Diego, CA, 1991, 151--176. 
\vskip1.5in 

\rightline{\sl July 9, 1992}

\end